\documentclass{jnmp}
\usepackage{amsmath}
\usepackage{amsthm}
\usepackage{amssymb}

\newcommand{\norm}[1]{\Vert #1 \Vert}  
\newcommand{\abs}[1]{\vert #1 \vert}  
\newcommand{\ov}{\overline}
\newcommand{\gr}{\mathrm{grad}}
\newcommand{\Ker}{\mathrm{Ker}\,}

\newcommand{\dif}{\mathrm{d}}
\newcommand{\tr}{\mathrm{trace}}
\newcommand{\CC}{\mathbb{C}}

\newcommand{\HH}{\mathcal{H}}
\newcommand{\VV}{\mathcal{V}}
\newcommand{\FF}{\mathcal{F}} 
 
\newcommand{\RR}{\mathbb{R}}

\newcommand{\Ss}{\mathbb{S}}

\JNMPnumberwithin{equation}{section}
\resetfootnoterule

\newtheorem*{acknowledgments}{Acknowledgments}

\newtheorem{pr}{Proposition}[section]
\newtheorem{co}{Corollary}[section]
\newtheorem{lm}{Lemma}[section]
\theoremstyle{definition}
\newtheorem{de}{Definition}[section]    
\newtheorem{re}{Remark}[section]
\newtheorem{ex}{Example}[section]

\begin{document}

\renewcommand{\evenhead}{R Slobodeanu}
\renewcommand{\oddhead}{Holomorphic maps and Faddeev-Hopf model}

\thispagestyle{empty}



\Name{A special class of holomorphic mappings and the Faddeev-Hopf model}

\label{firstpage}

\Author{Radu Slobodeanu}

\Address{Department of Theoretical Physics and Mathematics, University of Bucharest,\\
P.O. Box Mg-11, RO--077125 Bucharest-M\u agurele, Romania. \\
E-mail: radualexandru.slobodeanu@g.unibuc.ro}

\Date{July 18, 2013}


\begin{abstract}
\noindent
Pseudo horizontally weakly conformal maps \cite {lub} extend both holomorphic and (semi)conformal maps into an almost Hermitian manifold. We find critical points for the (generalized) Faddeev-Hopf model \cite {sve} in this larger class.
\end{abstract}

\section{Introduction}

Harmonic maps between Riemannian manifolds, $\varphi: (M,g) \to (N,h)$ are critical points for the Dirichlet energy functional (see e.g. \cite{ud}):

\begin{equation*}
\mathcal{E}(\varphi) =
\frac{1}{2} \int_M \norm{\dif \varphi}^2 v_g.
\end{equation*}
a generalization of the kinetic energy of Classical Mechanics.

Among the most studied classes of harmonic maps are those that are horizontally weakly conformal (HWC) (this is precisely the class of harmonic morphisms, cf. \cite{ud}). Another well known type of harmonic maps is given by holomorphic maps from a cosymplectic (semi-K\"ahler) manifold to a (1, 2)-symplectic (quasi-K\"ahler) one, according to a classical result of Lichnerowicz \cite{lih}. Moreover, holomorphic maps between almost K\"ahler manifolds are minimizers in their homotopy class for the Dirichlet energy.  

In this paper we turn our attention to the following  energy-type functional proposed in \cite{sve} as a generalization of the Faddeev-Hopf static Hamiltonian \cite{fad} from hadron physics:
\begin{equation} \label{efh}
\mathcal{E}_{\texttt{FH}}(\varphi) = \frac{1}{2} 
\int_M \left( \norm{\dif \varphi}^2 + 
\alpha \norm{\varphi^* \Omega}^2 \right) v_g,
\end{equation}
where $\alpha \geq 0$ is a coupling constant and $\Omega$ is the K\"ahler 2-form on $(N,J,h)$.

The main goal in what follows is to find critical maps for \eqref{efh} that belong to the two classes mentioned above (HWC and holomorphic). As the original model was concerned with the case $M=\mathbb{S}^3$, $N=\mathbb{S}^2$, we are interested in working with a broader notion of holomorphic map that allows also odd dimensional domains. Pseudo horizontally weakly conformal (PHWC) maps into almost Hermitian manifolds exhibit this extended holomorphicity and include also HWC mappings.

The paper is organized as follows. The next section surveys for further use some known facts about PHWC maps \cite{lub}, their harmonicity and the corresponding notion of pseudo horizontally homothetic (PHH) map. Besides we show that a wider perspective is possible by considering mappings from a Riemannian manifold to a $f$-structured one (rather than just almost complex) and by approaching their harmonicity also from the point of view of conformal geometry. We stress the fact that PHWC property answers the question when we can endow the domain with an almost $f$-structure compatible with the given metric, such that our map becomes holomorphic in an appropriate sense. Finally, Section 3 investigates the PHWC (stable) critical points for the strong coupling limit of the Faddeev-Hopf model using the variational formulae derived in \cite{sve}. Two of the  main applications are: 

\noindent (*) \ \textit{a holomorphic submersion between K\"ahler manifolds (or some appropriate analogue) is a critical point for the (full) Faddeev-Hopf functional if it has constant Dirichlet energy along horizontal curves}; 

\noindent (**) \ \textit{a K\"ahler targeted semiconformal (HWC) submersion is a critical point for the strong coupling Faddeev-Hopf functional if and only if it is 4-harmonic}.

Throughout this paper $M$, $N$ will be connected, $\mathcal{C}^{\infty}$ manifolds. All maps and other geometric objects considered will also be smooth.

\section{PHWC mappings} 

\subsection{Generalities about (almost) $f$-structures}
An almost $f$-structure on $M^m$ is a section $F$ of $\mathrm{End}(TM)$ such that $F^3+F = 0$. According to \cite{st}, the rank $k$, of $F$, is even and constant. If $k=m$, then $F$ is an \textit{almost complex structure} on $M$.

A pair $(g,F)$ is called a \textit{metric almost $f$-structure} if $g$ is a compatible metric, i.e. $g(FX,Y)+g(X,FY)=0$. In this case we have an \textit{orthogonal} decomposition:
$$
T^{\CC}M=T^0 M \oplus T^{(1,0)} M \oplus T^{(0,1)} M,
$$
where $T^0 M$, $T^{(1,0)} M$ and $T^{(0,1)} M$ are the eigenbundles of $F$ corresponding to 0, i and -i, respectively. We say that $\FF = T^0 M \oplus T^{(0,1)} M$ is \textbf{the complex distribution associated to} $F$, cf. \cite{lp}. 

\begin{de} (\cite{lp})
A map $\varphi : (M, F^M) \rightarrow (N, F^N)$, between manifolds endowed with almost $f$-structures, is called \textbf{holomorphic} if $\dif \varphi(\FF ^ M) \subseteq \FF ^ N$, where $\FF ^ M$ and $\FF ^ N$ are the complex distributions associated to $F^M$ and $F^N$, respectively.
\end{de}

In real terms, the above definition says that $\varphi$ is holomorphic if:
\begin{equation}
\dif \varphi(F^M X)-F^N \dif \varphi(X)\in \Ker F^N, \quad \forall X \in TM.
\end{equation}

In particular, when both $F^M$ and $F^N$ are almost complex structures, we find again the well-known definition of a holomorphic map: $\dif \varphi \circ J^M - J^N \circ \dif \varphi=0$.

\begin{de}(\cite{lp})
An almost $f$-structure $F$ is called \textbf{integrable} if $\FF$ is integrable (i.e. for any $X, Y \in \Gamma(\FF)$ we have $[X, Y ] \in \Gamma(\FF)$); an $f$-\emph{structure} is an integrable almost $f$-structure.
\end{de}

\begin{ex}
If an almost contact structure is \emph{normal} (\cite{ble}), then it is integrable in the sense of the above definition, cf. \cite{ia} (the converse is not necessarily true). For instance \emph{any} circle bundle over a complex manifold inherits on its total space, a normal almost contact structure. 
\end{ex}

\subsection{The general PHWC condition}

The following definition generalizes the corresponding one given in \cite {cen, lub} (when the codomain is endowed with an almost complex structure), respectively in  \cite{erd} (when the codomain is endowed with an almost $f$-structure):
\begin{de}
Let  $(M, g)$ be a Riemannian manifold and $(N, F, h)$ a manifold endowed with a metric almost $f$-structure. A map $\varphi: (M, g) \longrightarrow (N, F, h)$ is \textbf{pseudo horizontally weakly conformal} (PHWC) if:
\begin{equation} \label{PHWC}
F\circ [\dif \varphi \circ \dif \varphi^{t}, F]=0
\end{equation}
\noindent where $\dif \varphi^{t}$ stands for the \textit{adjoint map}, $\dif \varphi^{t}_x: T_{\varphi(x)}N \rightarrow T_x M$, characterized by: $g(X, \dif \varphi^{t}_x (E)) = h(\dif \varphi_x(X), E)$, $\forall X \in T_x M, E \in T_{\varphi(x)}N$. 
\end{de} 

The following result generalizes \cite[Remark 6]{lub}:
\begin{pr}
A constant rank PHWC map $\varphi: (M, g) \longrightarrow (N, F, h)$ induces a metric almost $f$-structure, $F^{\varphi}$, on $M$, with respect to which $\varphi$ becomes holomorphic.
\end{pr}

\begin{proof}
We note firstly that $\dif \varphi^{t}(T^{(1,0)}N)$ is a $g$-isotropic distribution on $M$.
\begin{equation*}
g(\dif \varphi^{t}(T^{(1,0)}N),\ \dif \varphi^{t}(T^{(1,0)}N))=h(T^{(1,0)}N,\ \dif \varphi \circ \dif \varphi^{t}(T^{(1,0)}N))=0,
\end{equation*}

as $\dif \varphi \circ \dif \varphi^{t}(T^{(1,0)}N)\subseteq T^{0}N \oplus T^{(1,0)}N$, cf. (\ref{PHWC}).
Therefore, according to \cite[Prop. 2.2]{ra}, there exists a unique metric $f$-structure $F^{\varphi}$ (of rank $2\mathrm{dim}[\dif \varphi^{t}(T^{(1,0)}N)]$) on $M$, such that $\dif \varphi^{t}(T^{(1,0)}N) = \Ker (F^{\varphi}-i)$. Denote by $T^{\CC}M = T^0 M \oplus T^{1,0} M \oplus T^{0,1} M$ the corresponding orthogonal splitting. It is easy to remark that: $\dif \varphi^{t}(T^{(0,1)}N) = \Ker (F^{\varphi}+i)$ and $\dif \varphi (T^0 M) \subseteq T^0 N$.

It is clear now that $\varphi$ is holomorphic, according to Definition 2.1, i.e. $\dif \varphi(T^0 M \oplus T^{(0,1)}M) \subseteq T^0 N \oplus T^{(0,1)}N$.
\end{proof}

\begin{re}
$(i)$ \ The PHWC condition does not depend on the metric on the codomain but only on the conformal class of the metric on the domain.

$(ii)$ \ $\mathrm{rank}F^{\varphi}=
\mathrm{rank}F+\mathrm{rank}\ \dif \varphi -
\mathrm{dim}N.$

In particular, if $\varphi$ is submersive ($\mathrm{rank} \ \dif \varphi = \mathrm{dim}N$), then $\mathrm{rank}F^{\varphi}=
\mathrm{rank}F$. If $F=J$ is an almost complex structure on $N$
($\mathrm{rank} F = \mathrm{dim}N$), then $\mathrm{rank}F^{\varphi}=\mathrm{rank} \ \dif \varphi$ (we have, moreover, $\Ker F^{\varphi}=\Ker \dif \varphi$). We call $F^{\varphi}$
the associated $f$-structure, cf. \cite{lub}.
\end{re}

\begin{ex}
Besides holomorphic maps between almost Hermitian manifolds, the following (already known) classes of mappings satisfy PHWC condition (\ref{PHWC}):

$(i)$ \emph{stable harmonic mappings to an irreducible Hermitian symmetric space of compact type}, cf. \cite {bur};

$(ii)$ $(\phi,J)$--\textit{holomorphic mappings from a metric almost contact manifold}  $(M, \phi, \xi, g)$ to an almost Hermitian manifold $(N,J,h)$, i.e. mappings such that $\dif \varphi \circ \phi = J \circ \dif \varphi$, cf. \cite{ian};

$(iii)$ \emph{semiconformal / horizontally weakly conformal (HWC) mappings}, i.e. mappings that satisfy $\dif \varphi \circ \dif \varphi^{t} = \lambda^{2} id$, cf. \cite{ud};

$(iv)$ \emph{contact - holomorphic mappings} between almost contact manifolds, cf. \cite{sl}.
\end{ex}

For the sake of simplicity, for the remainder of the article, we will make the assumption
that the target is endowed with an almost complex structure. Most of the results in the rest of this section can be easily extended to the general case of PHWC maps to a $f$--structured manifold.

\subsection{PHWC submersions}

Suppose now that $\varphi$ is submersive. The restriction of 
$\dif \varphi _x$ to the horizontal space 
$\HH _x = (\Ker \dif \varphi _x)^{\perp}$ maps that space isomorphically onto $T_{\varphi(x)}N$. Denote its inverse by $\widehat{\ }$ ; for any $Z\in T_{\varphi(x)}N$, the vector $\widehat{Z} \in \HH _x$ is called the \textit{horizontal lift of} $Z$ (this operation can be extended to local vector fields).

As $\varphi$ takes values in an almost Hermitian manifold $(N,J, h)$, we have an almost complex structure naturally induced on $\HH$ and an almost $f$-structure, $F$, that extends it: 
$F | _{\mathcal{H}}(X)=J_{\mathcal{H}}(X)=\widehat{J\dif \varphi(X)}$ and $F | _{\mathcal{V}}=0$, where $\VV=\Ker \dif \varphi$. 

In this case, the PHWC condition is equivalent to the compatibility of $J_{\HH}$ with the domain metric $g$ (i.e. $F$ is a metric almost $f$-structure). Indeed, starting with the remark that $\dif \varphi^t: TN \rightarrow \HH$ is this time an isomorphism, we obtain (for an arbitrary compatible metric $h$ on $N$):
\begin{equation*}
\begin{split}
g(J_{\mathcal{H}}X, \dif \varphi^t E)&=h(\dif \varphi(J_{\HH}X), E)=
h(J (\dif \varphi X), E)= -h(\dif \varphi X, JE)\\
&=-g(X, \dif \varphi^t JE)=-g(X, \dif \varphi^{-1}(\dif \varphi(\dif \varphi^t JE)))\\
&=-g(X, \dif \varphi^{-1}(J(\dif \varphi(\dif \varphi^t E)))=
-g(X, J_{\mathcal{H}}\dif \varphi^t E).
\end{split}
\end{equation*}

In this case, an alternative terminology for PHWC is \emph{horizontally holomorphic} cf. \cite{panti}.

\begin{re}
$(i)$ If $(N, J)$ is a \emph{complex} manifold (i.e. $J$ is integrable), then the almost $f$-structure induced by a PHWC submersion $\varphi: (M, g) \longrightarrow (N, J)$ is integrable, according to \cite[Remark 2.2]{pant}. Moreover,  \emph{all} $f$-structures (locally) appear in this way. 

\medskip
$(ii)$ Let $(z^{\alpha})_{\alpha=1,...,n}$ be local complex coordinates on the complex manifold $(N,J)$. The PHWC condition is equivalent to $g^{ij} \tfrac{\partial \varphi^{\alpha}}{\partial x ^{i}} \tfrac{\partial \varphi^{\beta}}{\partial x ^{j}}=0, \ \forall \alpha, \beta$, cf. \cite[Lemma 3]{lub}.

\medskip
$(iii)$
If the fibers of $\varphi$ are 1-dim. and $M$ is orientable, we have seen \cite{slob} that the associated $f$-structure is an almost contact metric structure, $\phi$. It is integrable but not necessarily normal. The supplementary condition that assures the normality is given in \cite[Theorem 4.1]{slob}.
\end{re}

\subsection{Harmonicity of PHWC maps}

A mapping is harmonic if and only if its tension field vanishes, cf. \cite{ud}. The tension field of a (constant rank) PHWC map into an almost Hermitian manifold $(N,J,h)$ is given by
\begin{equation} \label{ta}
\tau(\varphi)=J \mathrm{div}^{\varphi} J - \dif \varphi(F \mathrm{div} F), 
\end{equation}
where $F$ is the associated $f$-structure, $\mathrm{div}^{\varphi} J = \tr _g \varphi^{*}\nabla^N J$ and $\mathrm{div} F$ = trace $\nabla F$.

Therefore a PHWC mapping to a (1,2)-symplectic manifold (i.e. $(\nabla_{X}J)(Y)+(\nabla_{JX}J)(JY)=0$) will be harmonic if and only if $F \mathrm{div} F=0$, cf. \cite{mo}. 
If $F \mathrm{div} F=0$, then we shall call the (almost) $f$-structure, $F$,  \textit{cosymplectic}. 

For submersions, if we consider a local adapted frame $\{e_{i}, F e_{i}, e_{\alpha}\}$ (i.e. an orthonormal frame such that  $e_{\alpha} \in \Ker{F}, \forall \alpha$), then the above relation reads: 
\begin{equation} \label{tau}
\begin{split}
\tau(\varphi) =&  \sum_{i=1}^{n} J\left((\nabla^{\varphi}_{e_{i}}J)(\dif \varphi(e_{i}))+ (\nabla^{\varphi}_{Fe_{i}}J)(\dif \varphi(Fe_{i}))\right)\\
&- \dif \varphi\left(F \left[(\nabla_{e_{i}}F)(e_{i})+ (\nabla_{F e_{i}}F)(Fe_{i})\right]+(m-2n)\mu^{\mathcal{V}} \right),
\end{split}
\end{equation}
where $\mu^{\mathcal{V}}$ denotes the mean curvature of the fibres.

\subsection{Pseudo horizontally homothetic maps} \label{phhm}

Let $\varphi: (M, g) \longrightarrow (N, J, h)$ be a PHWC submersion with minimal fibres onto a (1,2)-symplectic manifold. If the associated $f$-structure $F$, on $M$, satisfies
\begin{equation}\label{div}
F\left((\nabla_{X}F)(X)+(\nabla_{FX}F)(FX)\right)=0, \quad \forall X \in \Ker(F^2+I),
\end{equation}
then $\varphi$ is harmonic and it has the following property: if $P$ is an complex submanifold of $N$, then $K=\varphi^{-1}(P) \subset M$ is a $f$-invariant minimal submanifold of $M$, cf. \cite[Prop. 7]{mo}. 

In \cite {aab}, the same properties were proved for \emph{pseudo horizontally homothetic} (PHH) submersions with minimal fibers. Recall that a PHWC map is PHH if it satisfies:
\begin{equation*}
[\dif \varphi \circ \nabla^{M}_{X} \circ \dif \varphi^{t}, J]=0, \quad \forall X \in \Ker(F^2 + I).
\end{equation*}
For submersions, this means that $J_{\mathcal{H}}$ is $\nabla^{\HH}$-parallel in horizontal directions (or, equivalently
$F((\nabla_{X}F)(Y))=0, \ \forall X,Y \in \Ker(F^2+I)$). In particular, a $(\phi,J)$-holomorphic submersion defined on a Sasakian manifold or a holomorphic submersion defined on a K\"ahler manifold are both PHH.

For horizontally weakly conformal maps both conditions \eqref{div} and PHH reduce to horizontal homothety (HH), i.e. the conformal factor $\lambda$ is constant in horizontal directions. Further properties and examples of PHH harmonic submersions can be found in \cite {aa, brz}.

\subsection{Conformal geometry viewpoint}

Let $(M,[g])$ be a conformal manifold ($[g]$ denotes a conformal class of Riemannian metrics). 

A \textbf{Weyl connection} $D$ on $(M,[g])$ is a torsion-free linear connection which preserves the conformal class $[g]$, cf. \cite{ga} (in this case we say that $D$ defines a \textit{Weyl structure} on $M$).
Preserving the conformal class means that for any $\tilde g \in [g]$, there exists a 1-form $\theta_{\tilde g}$ (called the \textit{Higgs field}) such that:
\begin{equation*}
D \tilde g = -2\theta_{\tilde g} \otimes \tilde g.
\end{equation*}
This formula is conformally invariant in the sense that,
if $\tilde g = e^{2f} g$, then $\theta_{\tilde g} = \theta_g - \dif f$.

Conversely, if one starts with a fixed Riemannian metric $g$ on $M$ (with $\nabla^g$ the Levi-Civita connection) and a fixed 1-form $\theta$ (with $\theta^{\sharp}$ the dual vector field with respect to $g$), then the connection
\begin{equation}\label{d}
D_{X}Y=\nabla^{g}_{X}Y + \theta(X)Y + \theta(Y)X - g(X,Y) \theta^\sharp,
\end{equation}
is a Weyl connection, preserving the conformal class of $g$. Clearly, $(g, \theta)$ and $(e^{2f} g , \theta - \dif f)$ define the same Weyl structure. If, moreover, $\theta$ is an exact 1-form, then $D$ is the Levi-Civita connection of some representant of the conformal class.

\medskip 

Suppose now $(M^m, g)$ endowed with a metric almost $f$-structure. Consider the Weyl structure (\ref{d}) constructed with an arbitrary 1-form on $M$. Then, with respect to an adapted orthonormal frame $\{ e_{i}, F e_{i}, e_{\alpha} \}$ (i.e. $e_{\alpha}$'s span $\mathrm{Ker} F$), we have:

\begin{equation*}
\begin{split}
F\mathrm{div}^{D}F &= \sum_{i, \alpha} F \left[(D_{e_{i}}F)(e_{i})+ (D_{F e_{i}}F)(Fe_{i}) + (D_{e_{\alpha}}F)(e_{\alpha}) \right]\\
&=F\mathrm{div}^{\nabla^g}F +(m-2)F^2 \theta^\sharp.
\end{split}
\end{equation*}

So if we take $\theta^\sharp=\frac{1}{m-2} F\mathrm{div}^{\nabla^g}F$, then $D$, the Weyl connection associated to the dual 1-form $\theta$, will have the property that:
\begin{equation*}
F\mathrm{div}^{D}F =0.
\end{equation*}

Notice that $D$ is not uniquely determined. Analogously to \cite[Def. 4.2]{lpa} for the almost complex case, we introduce the following: 

\begin{de}
A Weyl connection on $(M, [g], F)$ will be called \textbf{Weyl connection compatible with F} if $F\mathrm{div}^{D}F =0$.
\end{de}

Recall that we have an extended notion of harmonicity for maps defined on a conformal manifold \cite[Def. 2.1]{lpa}. Analogously to \cite[Prop. 4.5]{lpa}, we have:

\begin{pr}
A PHWC map $\varphi: (M,g) \rightarrow (N,J,h)$ to a $(1,2)$-symplectic manifold is harmonic with respect to some compatible Weyl connection on $(M, [g], F^{\varphi})$.
\end{pr}

\section{Critical points for the Faddeev-Hopf functional}

The strong coupling limit for the (generalized) \textit{Faddeev-Hopf model} involves the variational problem for the following energy-type functional, cf. \cite{sve}:
\begin{equation}\label{fhinfty}
\mathcal{E}_{\texttt{FH}}^{\infty}(\varphi) :=
\lim_{\alpha \to \infty}\alpha^{-1}\mathcal{E}_{\texttt{FH}}(\varphi)=
\norm{\varphi^* \Omega}_{L^2}^2
= \frac{1}{2} \int_M \langle \varphi^* \Omega, \ \varphi^* \Omega \rangle v_g,
\end{equation}
where $\varphi: (M, g) \rightarrow (N, J, h)$ are mappings defined on a compact, oriented Riemannian manifold, taking values in a K\"ahler manifold with the \textit{fundamental $2$-form} $\Omega=h(J\cdot, \cdot)$. According to \cite[Corollary 2.4]{sve}, such a mapping will be a critical point for this functional if and only if:
\begin{equation}\label{fh}
(\delta \varphi^* \Omega)^{\sharp}
\in \Ker \dif \varphi.
\end{equation}

We can remark that if the target is only a symplectic manifold, the above result is still true (so in the following the K\"ahler hypothesis can be weakened). 

The goal of this section is to identify PHWC submersions that satisfy (\ref{fh}). Unless otherwise stated, for the remainder of the article we assume that all mappings are K\"ahler targeted. 

Let us begin with an easy to check formula that will be useful in what follows:

\begin{lm}
Let $\varphi:(M,g) \to (N,h)$ be a mapping between Riemannian manifolds. Then:
\begin{equation}
\left( \nabla_X \varphi^* h \right)(Y, Z)= h(\nabla \dif \varphi(X, Y), \dif \varphi(Z)) + h(\dif \varphi(Y), \nabla \dif \varphi(X,Z)), \quad \forall X,Y,Z \in \Gamma(TM).
\end{equation}
\end{lm}

\begin{pr}\label{pr41}
Let $\varphi: (M, g) \rightarrow (N^{2n}, J, h)$ be a PHWC submersion from a compact, oriented Riemannian manifold. Then any two of the following statements imply the third:

\medskip
$(i)$ The $f$-structure, $F$, induced on $M$ is cosymplectic  (i.e. $F \mathrm{div}F=0$) 

\medskip
$(ii)$ $\varphi$ is a critical point of $\mathcal{E}_{\texttt{FH}}^{\infty}$ 
(i.e. it verifies the equation (\ref{fh}))

\medskip
$(iii)$ With respect to any (local) adapted orthonormal frame $\{E_{j}, F E_{j}, E_{\alpha}\}$ on $M$, the pullback of the codomain metric satisfies:
\begin{equation*}
\sum_{j=1}^{n} \left[ \left( \nabla_{E_j} \varphi^* h \right)(FE_j, Z)-\left( \nabla_{FE_j} \varphi^* h \right)(E_j, Z)\right]=0, \ \forall Z \in (\Ker \dif \varphi)^{\perp}.
\end{equation*}
\end{pr}

\begin{proof}
The PHWC condition assures us that: $\hat \Omega(X,Y):= \varphi^* \Omega (X, Y)=\varphi^* h (FX, Y)$. In order to compute its co-differential, we start with:
\begin{equation*}
\begin{split}
(\nabla_{X} \hat \Omega) (Y, Z) 
&= X[\varphi^* h(FY, Z)]-
\varphi^* h(F \nabla_{X} Y, Z)-\varphi^* h(FY, \nabla_X Z)\\
&= [h(\nabla^{\varphi}_{X} \dif \varphi(FY), \dif \varphi(Z))+
h(\dif \varphi(FY), \nabla^{\varphi}_{X} \dif \varphi(Z))]\\
&-\varphi^* h(F \nabla_{X} Y, Z)-\varphi^* h(FY, \nabla_X Z)\\
&=\varphi^* h((\nabla_{X} F)Y, Z)+
h(\nabla \dif \varphi (X, FY), \dif \varphi(Z))+
h(\dif \varphi(FY), \nabla \dif \varphi(X,Z))\\
&=\varphi^* h((\nabla_{X} F)Y, Z)+
\left( \nabla_X \varphi^* h \right)(FY, Z),
\quad \forall X,Y, Z \in \Ker (F^2 + I).
\end{split}
\end{equation*}

With respect to a (local) adapted orthonormal frame $\{E_{j}, F E_{j}, E_{\alpha}\}$ (i.e. $E_{\alpha} \in \Ker F$), for any horizontal vector field $Z$ on $M$, we have:
\begin{equation}
\begin{split}
-\delta \hat \Omega (Z)
&= \sum_{j,\alpha} \imath_{E_{j}}(\nabla_{E_{j}} \hat \Omega)(Z)+
\imath_{F E_{j}}(\nabla_{F E_{j}} \hat \Omega)(Z)+
\imath_{E_{\alpha}}(\nabla_{E_{\alpha}} \hat \Omega)(Z) \\
&=\varphi^* h(\mathrm{div}F, Z)+
\sum_{j} \left[ \left( \nabla_{E_j} \varphi^* h \right)(FE_j, Z)-\left( \nabla_{FE_j} \varphi^* h \right)(E_j, Z)\right].
\end{split}
\end{equation}

On the other hand, $\delta \hat \Omega (Z)=g((\delta \varphi^* \Omega)^{\sharp}, Z)$. As $\varphi$ is a submersion, $\Ker \dif \varphi = \Ker F$ and our conclusion easily follows.
\end{proof}

Note that any submersive harmonic PHWC map which satisfies $(iii)$ from Proposition \ref{pr41} will be critical point of the \textit{full} Faddeev-Hopf functional (\ref{efh}) (i.e. for every value of the coupling, $\alpha$, not just the infinite coupling limit).

\begin{co}
A PHWC submersion that satisfies
\begin{equation}\label{cond}
\nabla \dif \varphi (X, \ T^{(0,1)}M) \subseteq \varphi^{-1}T^{(1,0)}N, \quad \forall X \in \Gamma(\HH)
\end{equation}
is a critical point for $\mathcal{E}_{\texttt{FH}}^{\infty}$ if and only if it has minimal fibres.
\end{co}

\begin{proof}
The condition (\ref{cond}) assures us that 
$\nabla \dif \varphi (X, FY)=-J\nabla \dif \varphi (X, Y)$ which implies immediately the condition $(iii)$ from Proposition \ref{pr41}.

On the other hand, as $N$ is endowed with a K\"ahler structure,
from the following easy to check relation (true for all $X,Y \in \Gamma(\HH)$):
\begin{equation}\label{ndf}
0=\left( \nabla^{\varphi}_{X}J \right) \dif \varphi (Y) =
\dif \varphi \left( (\nabla_{X}F)Y \right)+\nabla \dif \varphi (X, FY)- 
J \nabla \dif \varphi (X,Y), 
\end{equation}
we can deduce that $\dif \varphi \left( (\nabla_{X}F)Y \right)+2\nabla \dif \varphi (X, FY)=0$. In particular, we also have: 

\noindent $\dif \varphi \left( (\nabla_{X}F)X + (\nabla_{FX}F)FX \right) = 0$ and therefore the condition $(i)$ from Proposition \ref{pr41} is also satisfied if and only if $\varphi$ has minimal fibres.
\end{proof}

Notice that, in order to have \eqref{cond}, $\nabla \dif \varphi$ must have no (1,1)-part (i.e. is $(1,1)$-geodesic, see \cite{ud}) and moreover $\nabla \dif \varphi(\ov Z, \ov W) \in \varphi^{-1}T^{(1,0)}N, \forall Z,W \in T^{(1,0)}M$. In the integrable case (e.g. for a holomorphic map between Hermitian manifolds) these conditions force the map to be totally geodesic (i.e. $\nabla \dif \varphi =0$).

As the strong coupling term of the Faddeev model comprises fourth power derivative terms, we expect some resemblance to the case of 4--energy. Recall (\cite{take}) that the $p$--\textit{energy} is defined as $\mathcal{E}_p(\varphi) = \frac{1}{p} \int_M \norm{\dif \varphi}^p v_g$ and its critical points are called $p$--\textit{harmonic maps}. If, in addition, a $p$--harmonic map is horizontally weakly conformal, then it pulls (local) $p$-harmonic functions on $N$ back to (local) $p$-harmonic functions on $M$ and it is called $p$-harmonic morphism, cf. \cite{lube}.

\begin{co}\label{32}
$(i)$ \ A PHH submersion is a critical point for the Faddeev-Hopf functional $\mathcal{E}_{\texttt{FH}}^{\infty}$ if and only if  $\mathrm{div}\varphi^* h | _{\HH} =0$.

\noindent $(ii)$ A PHH submersion with $\mathrm{grad}^{\HH}\norm{\dif \varphi}^2 =0$ is a critical point of the Faddeev-Hopf functional $\mathcal{E}_{\texttt{FH}}^{\infty}$ if and only if it has minimal fibres. In this case it is moreover harmonic and $4$-harmonic.
\end{co}

\begin{proof}
By a similar computation as for the proof of 
Proposition \ref{pr41}, using \eqref{ndf} we can check that:
\begin{equation*}
(\nabla_{X} \hat \Omega) (Y, Z) 
=-\left( \nabla_X \varphi^* h \right)(Y, FZ)
-\varphi^* h(Y, (\nabla_{X} F) Z),
\quad  \forall X,Y, Z \in \Ker (F^2 + I).
\end{equation*}

Therefore:
\begin{equation*}
\delta \hat \Omega (Z)
=\mathrm{div}\varphi^* h(FZ)+
\sum_{j} \left[ \varphi^* h(E_j, (\nabla_{E_j} F) Z)
+ \varphi^* h(FE_j, (\nabla_{FE_j} F) Z)\right].
\end{equation*}

As our PHH hypothesis (i.e. $F$ parallel in horizontal directions) assures the cancellation of every term in the above sum, the conclusion follows.

For the second statement, use the fact that a PHH submersion is harmonic if and only if it has minimal fibres and then take into account that a submersion is harmonic if and only if its stress-energy tensor is divergence free: $\mathrm{div}S_{\varphi} = (1/2) \dif (\norm {\dif \varphi}^2) - \mathrm{div}\varphi^{*}h =  0$.
\end{proof}

Notice that the hypothesis $\mathrm{grad}^{\HH}\norm{\dif \varphi}^2 =0$ is obviously satisfied by HH submersions, totally geodesic maps and eigenmaps between spheres. Moreover, it tells us that the submersion $\varphi$ must be $\infty$-\textit{harmonic}, a notion introduced in  \cite{ou, ouu}, where various examples are constructed and some classification results are given. Particularly, holomorphic $\infty$-harmonic maps $\mathbb{C}^n \to \mathbb{C}$ are a composition of an orthogonal projection followed by a homothety, cf. \cite{ou}. 

Recalling subsection \ref{phhm}, we see that Corollary \ref{32} applies to the following classes of PHH harmonic submersions.

\begin{ex}\label{BWfibr}
$(i.)$ A holomorphic submersion from a K\"ahler manifold (or a compact Vaisman manifold) onto a K\"ahler manifold, which has constant energy density (in horizontal directions), is a critical point of the full Faddeev-Hopf functional \eqref{efh}.

$(ii.)$ A Boothby-Wang fibration of a compact, regular Sasakian (or just K-contact) manifold over a K\"ahler (or just almost K\"ahler) manifold is a critical point of the full Faddeev-Hopf functional \eqref{efh} (on the total space we consider the metric $g=\varphi^* h + \eta \otimes \eta$ induced from the base, cf. \cite{ble}, so the fibration becomes a Riemannian submersion and therefore $\norm{\dif \varphi}$ is constant).
\end{ex}

For semiconformal (HWC) particular case, a more precise statement can be made, providing us with class of examples that extend \cite[Examples 3.1, 3.3]{sve}:

\begin{pr}
A semiconformal (HWC) submersion $\varphi: (M^m, g) \rightarrow (N^{2n}, J, h)$ with dilation $\lambda$ is a critical point for the Faddeev-Hopf functional $\mathcal{E}_{\texttt{FH}}^{\infty}$ if and only if:
\begin{equation}
(2n-4)\mathrm{grad}^{\HH}(\mathrm{ln} \lambda)+(m-2n)\mu^{\VV}=0.
\end{equation}
that is, if and only if it is $4$-harmonic, so a $4$-harmonic morphism.
\end{pr}

\begin{proof}
Semiconformal submersions are in particular PHWC, so we always have an induced metric almost $f$-structure, $F$, on $M$ and, in addition:
\begin{equation} \label{fla}
F\mathrm{div}F=(2n-2)\mathrm{grad}^{\HH}(\mathrm{ln} \lambda)+(m-2n)\mu^{\VV}.
\end{equation}
As $\varphi^* h \vert _ {\HH \times \HH} = \lambda^2 g \vert _ {\HH \times \HH}$, in this case we shall have, for any $Y$, $Z$ horizontal vector fields :
\begin{equation*}
(\nabla_{X} \hat \Omega) (Y, Z) 
= X(\lambda^2)g(FY, Z)+\lambda^2 g((\nabla_{X} F)Y, Z).
\end{equation*}
Then, taking the trace with respect to an adapted orthonormal frame, we obtain:
\begin{equation}
-\delta \hat \Omega(Z) = \lambda^2 g(F\mathrm{div}F -2\mathrm{grad}^{\HH}(\ln \lambda), FZ)
\end{equation}
Now taking into account also (\ref{fla}), the conclusion follows.
\end{proof}

\begin{co}\label{cohh}
$(i)$. A horizontally homothetic submersion is a stationary point for the \emph{full} Faddeev-Hopf functional \eqref{efh} if and only if it has minimal fibres. 

$(ii)$. A semiconformal submersion onto a four-manifold is a stationary point for the \emph{full} Faddeev-Hopf functional \eqref{efh}
if and only if it has minimal fibres.
\end{co}

\subsection{Stability}

According to \cite[Corollary 4.9]{sve}, the Hessian of a critical $\varphi$ for the energy $\mathcal{E}_{\texttt{FH}}^{\infty}$ is
\begin{equation*}
\mathrm{Hess}_{\varphi}(v,v)=\norm{\dif (\varphi^* \imath _v \Omega)}^{2}_{L^2}+\int _M \Omega(v, \nabla^{\varphi}_{Z_\varphi} v) \nu_g,
\quad \forall v \in \Gamma(\varphi^{-1}TN)
\end{equation*}
where $Z_\varphi=(\delta \varphi^* \Omega)^{\sharp}$. 

In particular, if $\delta \varphi^* \Omega (V)=0, \forall V \in \Ker \dif \varphi$, then $\varphi$ is (weakly) stable. But for a PHWC submersion we can check that 
$(\nabla_X \varphi^* \Omega)(Y,V)=-\varphi^* h(FY, \nabla_X V)$ and therefore
\begin{equation*}
-\delta \varphi^* \Omega (V)=\sum_{i=1}^{n}\lambda_{i}^{2}g([E_i, FE_i], V), \quad \forall V \in \Ker \dif \varphi,
\end{equation*}
where $\lambda_{i}^{2}$ are the nonzero eigenvalues of $\varphi^* h$ with respect to $g$ and $\{E_i, FE_i\}$ a frame of horizontal eigenvector fields (which exists around almost every point of $M$, cf. \cite{panti}).
We can easily conclude that a critical PHWC submersion for  $\mathcal{E}_{\texttt{FH}}^{\infty}$ is (weakly) stable if any of the following statements holds good:

\noindent $(a)$ \ the horizontal distribution $\HH$ is integrable;

\noindent $(b)$ \ the associated $f$-structure on $M$ satisfies
\begin{equation*}
(\nabla_{X}F)(X)+(\nabla_{FX}F)(FX)=0, \quad \forall X \in \Gamma(\HH).
\end{equation*}

\noindent Let us illustrate these situations (for the contact geometry background we refer to \cite{ble}).

\begin{pr}
Let $\varphi: (M^{2n+1}, \phi, \xi, \eta, g) \to (N^{2n},J,h)$ be a $(\phi, J)$-holomorphic, horizontally homothetic submersion from an almost contact metric manifold to a (almost) K\"ahler one. Then $\varphi$ is a stable critical point for $\mathcal{E}_{\texttt{FH}}^{\infty}$ in any of the following cases:

\noindent $(i)$ $(\phi, \xi, \eta, g)$ is a nearly cosymplectic structure;

\noindent $(ii)$ $(\phi, \xi, \eta, g)$ is a Kenmotsu structure.
\end{pr}

\begin{proof}
In both cases we can apply Corollary \ref{cohh}-$(i)$ to conclude that the map is critical.

Nearly cosymplectic structures are defined by $(\nabla_X \phi)X=0$, so the above condition $(b)$ is true, providing stability.
As Kenmotsu manifolds are locally warped products of an open interval with a K\"ahler manifold, $\Ker \eta$ (which must coincide with $\HH$) is integrable, so in this case the above condition $(a)$ is true.
\end{proof}


Recall that the Hopf map $\mathbb{S}^3 \to \mathbb{C}P^1$ is a stable critical point for $\mathcal{E}_{\texttt{FH}}^{\infty}$, according to \cite[Theorem 5.2]{sve}. As it is in particular a Boothby-Wang fibration, we are motivated to check the stability of this class of critical maps, considered in Example \ref{BWfibr}-$(ii)$.

Let $\varphi: (M^{2n+1}, \phi, \xi, \eta, g) \to (N^{2n},J,h)$ be a $(\phi, J)$-holomorphic, Riemannian submersion from a compact Sasakian manifold to a K\"ahler one. As $\xi$ must be in the kernel of $\dif \varphi$, we have
\begin{equation*}
Z_{\varphi}=
\delta \varphi^* \Omega (\xi)\xi=-2n \xi.
\end{equation*}

As $\varphi$ is submersive, for any $v \in \Gamma(\varphi^{-1}TN)$ there exists a local horizontal vector field $X_v$ on $M$ such that $v = \dif \varphi(X_v)$. But for all $X$ we have $\nabla^{\varphi}_{\xi} \dif \varphi(X) = \dif \varphi([\xi,X])$, so 
\begin{equation*}
\Omega(v, \nabla^{\varphi}_{Z_{\varphi}} v)=-2n\varphi^* h(\phi X_v, [\xi, X_v])=-2n g(\phi X_v, [\xi, X_v]).
\end{equation*}

Making explicit the term $\norm{\dif (\varphi^* \imath _v \Omega)}^{2}_{L^2}$ also, we get 
\begin{equation}\label{hesasa}
\begin{split}
\mathrm{Hess}_{\varphi}(v,v) &=
\int _M \left\{ \frac{1}{2}\sum_{\abs{I}\neq \abs{J}}\dif (\varphi^*\imath _v \Omega)(e_I, e_J)^{2} + (\mathrm{div}X_v)^2 + \norm{[\xi, X_v]}^2-2n g(\phi X_v, [\xi, X_v]) \right\} \nu_g , 
\end{split}
\end{equation}
where $\{ \xi, e_I \}_{I=1,...,n, \ov 1, ..., \ov n}$ is a local orthonormal frame on $M$ ($e_{\ov \imath}:=\phi e_i$).
We illustrate this situation with the following:

\begin{pr}
If $n \geq 2$, then the Hopf map $\mathbb{S}^{2n+1} \to \mathbb{C}P^n$ is an unstable critical point of $\mathcal{E}_{\texttt{FH}}^{\infty}$.
\end{pr}

\begin{proof}
Let $\varphi:\mathbb{S}^{2n+1} \to \mathbb{C}P^n$ be the Hopf map and  $(a_\alpha)_{\alpha=1,...,2n+2}$ be an orthonormal basis in $\RR^{2n+2}$. Define $f_\alpha : \Ss^{2n+1} \to \RR$, $f_\alpha(x)=\langle a_\alpha, x\rangle$ and $v_\alpha=\dif \varphi(\gr f_\alpha) \in \Gamma(\varphi^{-1}T\mathbb{C}P^n)$.
We have $(\gr f_\alpha)_ {x} = a_\alpha - f_\alpha(x) x$, $\abs{\gr f_\alpha}^2=1 - f_\alpha^2$ and
\begin{equation}\label{nabgr}
\nabla_X \gr f_\alpha=-f_\alpha X, \qquad (X \in \Gamma(T\Ss^{2n+1})),
\end{equation}
where $\nabla$ is the Levi-Civita connection on $\Ss^{2n+1}$.
Using the fact that $\Ss^{2n+1}$ is endowed with a natural Sasakian structure (so the above paragraph holds in this case) we compute
\begin{equation}\label{h1}
\Omega(v_\alpha, \nabla^{\varphi}_{Z_{\varphi}} v_\alpha)=-2n g(\phi \gr f_\alpha, [\xi, \gr f_\alpha])=-2n \left(\abs{\gr f_\alpha}^2 - \xi(f_\alpha)^2\right),
\end{equation}
where we have used \eqref{nabgr} and the identitiy $\phi X = - \nabla_X \xi$ (true on any Sasakian manifold).

Moreover $\dif (\varphi^* \imath _{v_\alpha} \Omega)(X,Y)=-2f_\alpha\dif \eta(X, Y)$ and $\dif (\varphi^* \imath _{v_\alpha} \Omega)(\xi, X)=g(\gr f_\alpha, X)$ for any $X, Y \perp \xi$, where, in addition, we used that $\nabla^{\varphi}_{\xi} \dif \varphi(X) = \dif \varphi([\xi,X])$ and $\nabla^{\varphi}_{X} \dif \varphi(Y) = \dif \varphi(\nabla_X Y)$, $\varphi$ being a Riemannian submersion. Therefore
\begin{equation}\label{h2}
\abs{\dif (\varphi^* \imath _{v_\alpha} \Omega)}^{2}=\abs{\gr f_\alpha}^2 - \xi(f_\alpha)^2 + 4n f_\alpha^2.
\end{equation}

From Equations \eqref{h1} and \eqref{h2}, by taking the sum over $\alpha$, we obtain
\begin{equation}
\sum_\alpha \mathrm{Hess}_{\varphi}(v_\alpha, v_\alpha)=2n(3-2n)\mathrm{Vol}(\Ss^{2n+1}),
\end{equation}
and the conclusion follows.
\end{proof}

Notice that if $n \geq 2$, then the Hopf map is also unstable as 4-harmonic map, according to the general result in \cite{leu}.

\begin{acknowledgments}
I thank Eric Loubeau, Liviu Ornea and Radu Pantilie for useful remarks. This work was supported by the CEx grant no.\ 2-CEx 06-11-22/25.07.2006.
\end{acknowledgments}

\label{lastpage}


\begin{thebibliography}{99}

\bibitem {aa} \textsc{Aprodu M. A.} and \textsc{Aprodu M.}, Implicitly defined harmonic PHH submersions, \textit{Manuscripta Math.} {\bf 100} (1999), 103--121.

\bibitem {aab} \textsc{Aprodu M. A., Aprodu M.} and \textsc{Br\^{\i}nz\u anescu V.}, A class of harmonic submersions and minimal submanifolds, \textit{Int. J. Math.} {\bf 11}(9) (2000), 1177--1191. 

\bibitem {ud} \textsc{Baird P.} and \textsc{Wood J.C.}, Harmonic Morphisms Between Riemannian Manifolds, Clarendon Press - Oxford, 2003.


\bibitem {ble} \textsc{Blair D.E.}, Riemannian Geometry of Contact and Symplectic Manifolds, Birkhauser Boston, Progress in Mathematics, vol.203, 2002.

\bibitem {brz} \textsc{Br\^{\i}nz\u anescu V.}, Pseudo-harmonic morphisms; applications and examples, \textit{An. Univ. Timi\c soara Ser. Mat.-Inform.} \textbf{39} (2001), Special Issue: Mathematics, 111--121. 

\bibitem {sl} \textsc{Br\^{\i}nz\u anescu V.} and \textsc{Slobodeanu R.}, 
Holomorphicity and Walczak formula on Sasakian manifolds, \textit{J. Geom. Phys.} {\bf 57} (2006), 193--207.

\bibitem {bur} \textsc{Burns D., Burstall F., De Bartolomeis P.} and \textsc{Rawnsley J.}, Stability of harmonic maps of K\"ahler manifolds, \textit{J. Diff. Geom.}  {\bf 30}  (1989), 579--594.

\bibitem {cen} \textsc{Chen J.}, Structures of certain harmonic maps into K\"ahler manifolds, \textit{Int. J. Math.} {\bf 8} (1997), 573--581.

\bibitem {leu} \textsc{Cheung, L. F.} and \textsc{Leung P. F.}, Some results on stable $p$-harmonic maps, \textit{Glasgow Math. J.} \textbf{36} (1994), 77--80.

\bibitem {erd} \textsc{Erdem S.}, $\varphi$-pseudo harmonic morphisms, some subclasses and their liftings to tangent bundles, \textit{Houston J. Math.} \textbf{30} (2004), 1009--1038.

\bibitem {fad} \textsc{Faddeev L. D.} and \textsc{Niemi A. J.}, Stable knot-like structures in classical field theory, \textit{Nature} \textbf{387} (1997), 58--61.

\bibitem {ga} \textsc{Gauduchon P.}, Structures de Weyl et th\'eor\`emes d'annulation sur une vari\' et\' e conforme autoduale, \textit{Ann. Sc. Norm. Sup. Pisa} \textbf{XVIII 4} (1991), 563--629.

\bibitem {ia} \textsc{Ianu\c s S.}, Sulle variet\`a di Cauchy-Riemann, \textit{Rend. Dell'Academia de Scienze di Napoli}, \textbf{39} (1972), 191--195.

\bibitem {ian} \textsc{Ianu\c s, S.} and \textsc{Pastore, A.M.}, Harmonic maps on contact metric manifolds, \textit{Ann. Math. Blaise Pascal} {\bf 2} (1995),
 43--53.

\bibitem {lih} \textsc{Lichnerowicz A.}, Applications Harmoniques et Vari\'et\'es K\"ahleriennes, \textit{Sympos. Math.} {\bf III} (1970), 341--402. 

\bibitem {lub} \textsc{Loubeau E.}, Pseudo harmonic morphisms, \textit{Int. J. of Math.} {\bf 8} (1997), 943--957.

\bibitem {lube} \textsc{Loubeau E.}, On $p$-harmonic morphisms, \textit{Diff. Geom. and Appl.} {\bf 12} (2000), 219--229.

\bibitem {mo} \textsc{Loubeau E.} and \textsc{Mo X.}, The geometry of pseudo harmonic morphisms, \textit{Beitr\"age Algebra Geom.}  {\bf 45}(1) (2004), 87--102.

\bibitem {lpa} \textsc{Loubeau E.} and \textsc{Pantilie R.}, Harmonic morphisms between Weyl spaces and twistorial maps, \textit{Comm. Anal. Geom.} \textbf{14} (2006), 847--881.

\bibitem {lp} \textsc{Loubeau E.} and \textsc{Pantilie R.}, Harmonic morphisms between Weyl spaces and twistorial maps II, math.DG/0610676.

\bibitem {tond} \textsc{Nishikawa S.} and \textsc{Tondeur P.}, Transversal infinitesimal automorphisms for harmonic K\"ahler foliations, \textit{T\^ohoku Math. J.} {\bf 40} (1988), 599--611.

\bibitem {ou} \textsc{Ou Y.-L.} and \textsc{Wang, Z.-P.}, Some classifications of $\infty$-Harmonic maps between Riemannian manifolds, math.DG/0710.5752.

\bibitem {ouu} \textsc{Ou Y.-L., Troutman, T.} and \textsc{Wilhelm F.}, Infinity-harmonic maps and morphisms, math.DG/0810.0975.

\bibitem {pant} \textsc{Pantilie R.}, On a class of twistorial maps, \textit{Diff. Geom. Appl.} \textbf{26} (2008), 366--376.

\bibitem {panti} \textsc{Pantilie R.} and \textsc{Wood J.C.}, Harmonic morphisms with one-dimensional fibres on Einstein manifolds, \textit{Trans. Amer. Math. Soc.} \textbf{354} (2002), 4229--4243.

\bibitem {ra} \textsc{Rawnsley J. H.}, $f$-structures, $f$-twistor spaces and harmonic maps, in \textit{Geometry seminar "Luigi Bianchi" II}, Lecture Notes in Math., \textbf{1164}, Springer, 1985, 85--159.

\bibitem {slob} \textsc{Slobodeanu R.}, Pseudo-harmonic morphisms with low dimensional fibers, \textit{Rend. Circ. Mat. Palermo}, \textbf{II(LV)} (2006), 5--20.

\bibitem {sve} \textsc{Speight J. M.} and \textsc{Svensson M.}, On the Strong Coupling Limit of the Faddeev-Hopf Model, \textit{Commun. Math. Phys.} \textbf{272} (2007), 751--773.

\bibitem {st} \textsc{Stong R.E.}, The rank of an \textit{f}-structure, \textit{Kodai Math. Sem. Rep.} \textbf{29} (1977), 207--209.

\bibitem{take} \textsc{Takeuchi H.}, Stability and Liouville theorems of $p$-harmonic maps, \textit{Jap. J. Math.}, New Ser. \textbf{17} (1991), 317--332.
\end{thebibliography}
\end{document}